\documentclass[a4paper]{article}
\usepackage{amsthm}
\usepackage{fleqn}
\usepackage{array}
\usepackage[latin1]{inputenc}
\usepackage{amssymb}
\usepackage{amsfonts}
\usepackage{amsmath}

\begin{document}

\newtheorem{thm}{Theorem}[section]
\newtheorem{prop}{Proposition}[section]
\newtheorem{lem}{Lemma}[section]
\newtheorem{cor}{Corollary}[section]

\numberwithin{equation}{section}

\begin{center}
{\Large Rotational surfaces in a $3$-dimensional normed space}

\vspace{5mm}

Makoto SAKAKI
\end{center}

\vspace{3mm}

{\bf Abstract.} We study rotational surfaces with constant Minkowski Gaussian curvature and rotational surfaces with constant Minkowski mean curvature in a $3$-dimensional normed space with rotationally symmetric norm. We have a generalization of the catenoid, pseudo-sphere and Delaunay surfaces. 

\vspace{3mm}


\vspace{3mm}


\section {Introduction}

It is interesting to generalize differential geometry of curves and surfaces in Euclidean spaces to that in normed spaces, or generally, in gauge spaces (cf. \cite{BMSa1}, \cite{BMSa2}, \cite{BMSh}, \cite{BMT1}, \cite{BMT2}, \cite{BMT3}, \cite{BMT4}, \cite{Bu}, \cite{G}), where how to compensate for the lack of the notion of angle is the problem, and the notion of Birkhoff orthogonality plays an important role. For surfaces in $3$-dimensional normed spaces, the notions of Birkhoff-Gauss map, Minkowski Gaussian curvature and Minkowski mean curvature are particularly important (cf. \cite{BMT1}, \cite{BMT2}, \cite{BMT3}, \cite{BMT4}). 

In this paper, we study rotational surfaces in a $3$-dimensional normed space with rotationally symmetric norm, in particular, rotational surfaces with constant Minkowski Gaussian curvature and rotational surfaces with constant Minkowski mean curvature.

This paper is organized as follows. In Section 2, following \cite{BMT2}, we recall some basic facts on surfaces in $3$-dimensional normed spaces. In Section 3, we give a basic computation for rotational surfaces in a $3$-dimensional normed space with rotationally symmetric norm. In Section 4, we consider rotational minimal surfaces in the $3$-dimensional normed space. In Section 5, we discuss rotational surfaces with non-zero constant Minkowski Gaussian curvature in the $3$-dimensional normed space. In Section 6, we study rotational surfaces with non-zero constant Minkowski mean curvature in the $3$-dimensional normed space, which can be seen as a generalization of the Delaunay surfaces (\cite{D}).

\section{Preliminaries}

In this section, following \cite{BMT2}, we recall some basic facts on surfaces in $3$-dimensional normed spaces. 

Let $({\mathbb R}^3, \|\cdot\|)$ be a $3$-dimensional normed space whose unit ball $B$ and unit sphere $S$ are defined by
\[B = \{x \in {\mathbb R}^3; \| x\| \leq 1\}, \ \ \ \ S = \{x \in {\mathbb R}^3; \| x\| = 1\}.\]
In the following, we assume that $S$ is smooth and strictly convex, that is, $S$ is a smooth surface and $S$ contains no line segment.

\vspace{3mm}

Remark. We do not assume that $S$ has positive Euclidean Gaussian curvature as in \cite{BMT2}, because we treat the case where $S$ has points with zero Euclidean Gaussian curvature.

\vspace{3mm}

Let $v$ be a non-zero vector in ${\mathbb R}^3$ and $\Pi$ be a plane in ${\mathbb R}^3$. We say that $v$ is Birkhoff orthogonal to $\Pi$ (denoted by $v \dashv_{B} \Pi$) if the tangent plane of $S$ at $v/\|v\|$ is $\Pi$.

Let $M$ be a surface immersed in $({\mathbb R}^3, \|\cdot\|)$. Let $T_{p}M$ be the tangent plane of $M$ at $p \in M$. There exists a vector $\eta(p) \in S$ such that $\eta(p) \dashv_{B} T_{p}M$, which gives a local smooth map $\eta: U \subset M \rightarrow S$ called the Birkhoff-Gauss map. It can be global if and only if $M$ is orientable. We define the Minkowski Gaussian curvature $K$ and the Minkowski mean curvature $H$ of $M$ at $p$ by
\[K = \mbox{det}(d\eta_p), \ \ \ \ H = \frac{1}{2}\mbox{trace}(d\eta_p).\]
We say that $M$ is flat if $K = 0$ identically, and minimal if $H = 0$ identically. 

A surface which is homothetic to the unit sphere $S$ is called a Minkowski sphere. A Minkowski sphere has positive constant Minkowski Gaussian curvature and non-zero constant Minkowski mean curvature.

\section{Rotational surfaces}

Let
\[S = \{(x_1, x_2, x_3) \in {\mathbb R}^3 | (x_1^2+x_2^2)^m+x_3^{2m} = 1\} \]
where $m$ is a positive integer. It is given by rotating $x_1^{2m}+x_3^{2m} = 1$ around $x_3$-axis. Then there exists a norm $\|\cdot\|$ on ${\mathbb R}^3$ whose unit sphere is the above $S$. Set
\[\Phi(x_1, x_2, x_3): = (x_1^2+x_2^2)^m+x_3^{2m}. \]
Throughout this paper, we consider this $3$-dimensional normed space $({\mathbb R}^3, \|\cdot\|)$. The case where $m = 1$ is the Euclidean case. We assume that $m \geq 2$ in the following. 

Let $M$ be a surface in $({\mathbb R}^3, \|\cdot\|)$ which is rotational around $x_3$-axis, and is parametrized by
\[f(u, v) = (\alpha(u)\cos{v}, \alpha(u)\sin{v}, \beta(u)) \]
where $\alpha > 0$, $\alpha' \neq 0$ and $\beta' \neq 0$. Then
\[f_u = (\alpha'\cos{v}, \alpha'\sin{v}, \beta'), \ \ \ \ f_v = (-\alpha\sin{v}, \alpha\cos{v}, 0). \]
The Birkhoff-Gauss map $\eta$ is characterized by the condition
\[\mbox{grad}(\Phi)|_{\eta} = \mu f_u\times f_v, \]
where $\mu$ is a positive function and $\times$ denotes the standard cross product in ${\mathbb R}^3$. Then we can get
\[\eta = A^{-\frac{1}{2m}} \left( -(\beta')^{\frac{1}{2m-1}}\cos{v}, -(\beta')^{\frac{1}{2m-1}}\sin{v}, (\alpha')^{\frac{1}{2m-1}} \right) \]
where
\[A:= \left( \alpha' \right)^{\frac{2m}{2m-1}}+\left( \beta' \right)^{\frac{2m}{2m-1}}. \]

We can compute that
\begin{eqnarray}
\eta_u = -\frac{1}{2m-1} A^{-\frac{2m+1}{2m}} (\alpha')^{-\frac{2m-2}{2m-1}} (\beta')^{-\frac{2m-2}{2m-1}} (\alpha'\beta''-\alpha''\beta') f_u 
\end{eqnarray}
and
\begin{eqnarray}
\eta_v = -\frac{1}{\alpha} A^{-\frac{1}{2m}}(\beta')^{\frac{1}{2m-1}} f_v. 
\end{eqnarray}
Thus we have
\begin{eqnarray}
K = \frac{1}{(2m-1)\alpha} A^{-\frac{m+1}{m}} (\alpha')^{-\frac{2m-2}{2m-1}} (\beta')^{-\frac{2m-3}{2m-1}} (\alpha'\beta''-\alpha''\beta') 
\end{eqnarray}
and
\begin{eqnarray}
H = -\frac{1}{2(2m-1)\alpha} A^{-\frac{2m+1}{2m}} (\beta')^{-\frac{2m-2}{2m-1}} \nonumber
\end{eqnarray}
\begin{eqnarray}
\hspace{1.5cm} \times \left\{ \alpha (\alpha')^{-\frac{2m-2}{2m-1}} (\alpha'\beta''-\alpha''\beta') +(2m-1)A \beta' \right \}. 
\end{eqnarray}
Letting $\beta(u) = u$, we have
\begin{eqnarray}
\eta_u = \frac{1}{2m-1} \left( \left( \alpha' \right)^{\frac{2m}{2m-1}}+1 \right)^{-\frac{2m+1}{2m}} (\alpha')^{-\frac{2m-2}{2m-1}} \alpha'' f_u, 
\end{eqnarray}
\begin{eqnarray}
\eta_v = -\frac{1}{\alpha} \left( \left( \alpha' \right)^{\frac{2m}{2m-1}}+1 \right)^{-\frac{1}{2m}} f_v, 
\end{eqnarray}
\begin{eqnarray}
K = -\frac{1}{(2m-1)\alpha} \left( \left( \alpha' \right)^{\frac{2m}{2m-1}}+1 \right)^{-\frac{m+1}{m}} (\alpha')^{-\frac{2m-2}{2m-1}} \alpha'', 
\end{eqnarray}
and
\begin{eqnarray}
H = \frac{1}{2(2m-1)\alpha} \left( \left( \alpha' \right)^{\frac{2m}{2m-1}}+1 \right)^{-\frac{2m+1}{2m}} \nonumber
\end{eqnarray}
\begin{eqnarray}
\hspace{1.5cm} \times \left\{ \alpha (\alpha')^{-\frac{2m-2}{2m-1}} \alpha'' -(2m-1) \left( \left( \alpha' \right)^{\frac{2m}{2m-1}}+1 \right) \right\}. 
\end{eqnarray}
By (3.7), we see that $K = 0$ if and only if $\alpha'' = 0$. So we have the following.

\begin{prop}
A rotational surface in $({\mathbb R}^3, \|\cdot\|)$ parametrized by
\[f(u, v) = (\alpha(u)\cos{v}, \alpha(u)\sin{v}, u) \]
where $\alpha > 0$, $\alpha' \neq 0$ is flat if and only if it is a circular cone. 
\end{prop}

\section{Rotational minimal surfaces}

Let $({\mathbb R}^3, \|\cdot\|)$ be the $3$-dimensional normed space as in Section 3. Let $M$ be a rotational surface in $({\mathbb R}^3, \|\cdot\|)$ parametrized by
\[f(u, v) = (\alpha(u)\cos{v}, \alpha(u)\sin{v}, u) \]
where $\alpha > 0$ and $\alpha' \neq 0$. 

By (3.8), the rotational surface $M$ is minimal if and only if
\begin{eqnarray}
\alpha (\alpha')^{-\frac{2m-2}{2m-1}} \alpha'' -(2m-1) \left( \left( \alpha' \right)^{\frac{2m}{2m-1}}+1 \right) = 0. 
\end{eqnarray}
From the equation (4.1), we have
\[\frac{2m-1}{\alpha} = \frac{\alpha''}{(\alpha')^2+(\alpha')^{\frac{2m-2}{2m-1}}}.\]
Multiplying by $2\alpha'$ we have
\[2(2m-1)\frac{\alpha'}{\alpha} = \frac{((\alpha')^2)'}{(\alpha')^2+(\alpha')^{\frac{2m-2}{2m-1}}}, \]
and
\[2(2m-1)\log{\alpha} = \int \frac{((\alpha')^2)'}{(\alpha')^2+(\alpha')^{\frac{2m-2}{2m-1}}} du. \]
Letting
\[(\alpha')^{\frac{2}{2m-1}} =:Z \]
for the right-hand side, we have
\[2\log{\alpha} = \int \frac{Z^{m-1}}{Z^{m}+1}dZ = \frac{1}{m}\log{\left( Z^{m}+1 \right)}+c_1 \]
\[= \frac{1}{m}\log{\left( (\alpha')^{\frac{2m}{2m-1}}+1 \right)}+c_1 \]
for a constant $c_1$. Then
\[\frac{d\alpha}{du} = \pm \frac{1}{c_2^{2m-1}} \left( \alpha^{2m} -c_2^{2m} \right)^{\frac{2m-1}{2m}} \]
for a positive constant $c_2$, and
\[u(\alpha) = \pm \int_{c_2}^{\alpha} \frac{c_2^{2m-1}}{\left( \rho^{2m} -c_2^{2m} \right)^{\frac{2m-1}{2m}}} d\rho+c_3 \]
for a constant $c_3$, where $\alpha > c_2$. Here we note that since
\[0 < \frac{2m-1}{2m} < 1, \]
the above integral converges and
\[\lim_{\alpha\rightarrow c_2} \int_{c_2}^{\alpha} \frac{c_2^{2m-1}}{\left(\rho^{2m} - c_2^{2m} \right)^{\frac{2m-1}{2m}}} d\rho = 0. \]
Then we have the following. 

\begin{thm}
A rotational surface in $({\mathbb R}^3, \|\cdot\|)$ given by
\[\bar{f}(\alpha, v) = (\alpha \cos v, \alpha \sin v, u(\alpha)) \]
where $\alpha > 0$ is minimal if and only if
\[u(\alpha) = \pm \int_{c_2}^{\alpha} \frac{c_2^{2m-1}}{\left( \rho^{2m} -c_2^{2m} \right)^{\frac{2m-1}{2m}}} d\rho+c_3 \]
for constants $c_2 > 0$ and $c_3$, where $\alpha > c_2$. 
\end{thm}

Now, let us set
\[u_{\pm}(\alpha) := \pm \int_{c_2}^{\alpha} \frac{c_2^{2m-1}}{\left( \rho^{2m} -c_2^{2m} \right)^{\frac{2m-1}{2m}}} d\rho+c_3 \]
for constants $c_2 > 0$ and $c_3$, where $\alpha > c_2$, and we consider the behavior of the graph of $u_{\pm}(\alpha)$. Since $m \geq 2$, 
\[\lim_{\alpha\rightarrow \infty} \int_{c_2}^{\alpha} \frac{c_2^{2m-1}}{\left(\rho^{2m} - c_2^{2m} \right)^{\frac{2m-1}{2m}}} d\rho = d_1 \]
for some positive value $d_1$. So we have
\[\lim_{\alpha\rightarrow c_2} u_{\pm}(\alpha) = c_3, \ \ \ \ \lim_{\alpha\rightarrow \infty} u_{+}(\alpha) = c_3+d_1, \ \ \ \ \lim_{\alpha\rightarrow \infty} u_{-}(\alpha) = c_3-d_1. \]
The function $u_{+}(\alpha)$ is an increasing function and
\[\lim_{\alpha\rightarrow c_2} u_{+}'(\alpha) = \infty. \]
Similarly, $u_{-}(\alpha)$ is a decreasing function and
\[\lim_{\alpha\rightarrow c_2} u_{-}'(\alpha) = -\infty. \]
Let $\alpha_{+}(u)$ be the inverse function of $u_{+}(\alpha)$. It is an increasing function on $(c_3, c_3+d_1)$ and
\[\lim_{u\rightarrow c_3}\alpha_{+}(u) = c_2, \ \ \ \ \lim_{u\rightarrow c_3+d_1}\alpha_{+}(u) = \infty, \ \ \ \ \lim_{u\rightarrow c_3}\alpha_{+}'(u) = 0. \]
Let $\alpha_{-}(u)$ be the inverse function of $u_{-}(\alpha)$. It is a decreasing function on $(c_3-d_1, c_3)$ and
\[\lim_{u\rightarrow c_3}\alpha_{-}(u) = c_2, \ \ \ \ \lim_{u\rightarrow c_3-d_1}\alpha_{-}(u) = \infty, \ \ \ \ \lim_{u\rightarrow c_3}\alpha_{-}'(u) = 0. \]

We define a function $\hat{\alpha}(u)$ on $(c_3-d_1, c_3+d_1)$ by
\[\hat{\alpha}(u) = \left\{ \begin{array}{cc} 
\alpha_{+}(u), & c_3 < u < c_3+d_1 \\ \\
\alpha_{-}(u), & c_3-d_1 < u < c_3 \\ \\
c_2, & u = c_3 
\end{array}. \right. \]
Then $\hat{\alpha}(u)$ is a $C^1$-function on $(c_3-d_1, c_3+d_1)$ such that
\[\hat{\alpha}'(u) = \left\{ \begin{array}{cc} 
\alpha_{+}'(u), & c_3 < u < c_3+d_1 \\ \\
\alpha_{-}'(u), & c_3-d_1 < u < c_3 \\ \\
0, & u = c_3 
\end{array}. \right. \]
For $u \in (c_3-d_1, c_3) \cup (c_3, c_3+d_1)$, $\hat{\alpha}(u)$ satisfies the equation (4.1). Then, noting that $m \geq 2$, we can see that
\[\lim_{u\rightarrow c_3} (\hat{\alpha}'(u))^{-\frac{2m-2}{2m-1}} \hat{\alpha}''(u) = \frac{2m-1}{c_2} \]
and
\[\lim_{u\rightarrow c_3} \hat{\alpha}''(u) = 0. \]
Thus the function $\hat{\alpha}(u)$ is a $C^2$-function on $(c_3-d_1, c_3+d_1)$. Also by (3.5) and (3.6), we find that the Birkhoff-Gauss map can be $C^1$-extended for $u \in (c_3-d_1, c_3+d_1)$. 

Therefore, we have the following. 

\begin{thm}
Under the notation above, the rotational surface in $({\mathbb R}^3, \|\cdot\|)$ parametrized by
\[\hat{f}(u, v) = (\hat{\alpha}(u)\cos{v}, \hat{\alpha}(u)\sin{v}, u), \ \ \ \ (u, v) \in (c_3-d_1, c_3+d_1)\times [0, 2\pi] \]
is minimal. 
\end{thm}

Remark. The above surface can be seen as a generalization of the catenoid in the Euclidean $3$-space. But we should note that the range of $u$ is bounded.

\section{Non-zero constant Gaussian curvature}

Let $({\mathbb R}^3, \|\cdot\|)$ be the $3$-dimensional normed space as in Section 3. Let $M$ be a rotational surface in $({\mathbb R}^3, \|\cdot\|)$ parametrized by
\[f(u, v) = (\alpha(u)\cos{v}, \alpha(u)\sin{v}, u) \]
where $\alpha > 0$ and $\alpha' \neq 0$. 

By (3.7), if $K$ is a non-zero constant, then
\begin{eqnarray}
-\frac{1}{2m-1}\left( \left( \alpha' \right)^{\frac{2m}{2m-1}}+1 \right)^{-\frac{m+1}{m}} (\alpha')^{-\frac{2m-2}{2m-1}} \alpha'' = K\alpha. 
\end{eqnarray}
Multiplying by $2\alpha'$ we have
\[-\frac{2}{2m-1}\left( \left( \alpha' \right)^{\frac{2m}{2m-1}}+1 \right)^{-\frac{m+1}{m}} (\alpha')^{\frac{1}{2m-1}} \alpha'' = K(\alpha^2)'. \]
Integrating it we have
\[\left( \left( \alpha' \right)^{\frac{2m}{2m-1}}+1 \right)^{-\frac{1}{m}} = K\alpha^2+c_1 \ \ (> 0) \]
for a constant $c_1$. Then
\[\frac{d\alpha}{du} = \pm \frac{\left\{ 1-(K\alpha^2+c_1)^m \right\}^{\frac{2m-1}{2m}}}{(K\alpha^2+c_1)^{\frac{2m-1}{2}}}, \]
and we get the following. 

\begin{thm}
A rotational surface in $({\mathbb R}^3, \|\cdot\|)$ given by
\[\bar{f}(\alpha, v) = (\alpha \cos v, \alpha \sin v, u(\alpha)) \]
where $\alpha > 0$ has non-zero constant Minkowski Gaussian curvature $K$ if and only if
\[u(\alpha) = \pm \int \frac{(K\alpha^2+c_1)^{\frac{2m-1}{2}}} {\left\{ 1-(K\alpha^2+c_1)^m \right\}^{\frac{2m-1}{2m}}} d\alpha \]
for a constant $c_1$. 
\end{thm}

Now, let us set
\[u_{\pm}(\alpha) := \pm \int \frac{(K\alpha^2+c_1)^{\frac{2m-1}{2}}} {\left\{ 1-(K\alpha^2+c_1)^m \right\}^{\frac{2m-1}{2m}}} d\alpha \]
for a constant $c_1$, and we discuss the behavior of the graph of $u_{\pm}(\alpha)$. It suffices to consider the case where $K = 1$ or $K = -1$. 

\vspace{2mm}

(i) The case $K = 1$. We have $c_1 < 1$ and
\[u_{\pm}(\alpha) = \pm \int \frac{(\alpha^2+c_1)^{\frac{2m-1}{2}}} {\left\{ 1-(\alpha^2+c_1)^m \right\}^{\frac{2m-1}{2m}}} d\alpha. \]

(i-1) The case $c_1 = 0$. In this case we have
\[u_{\pm}(\alpha) = \pm \int \frac{\alpha^{2m-1}} {\left( 1-\alpha^{2m} \right)^{\frac{2m-1}{2m}} } d\alpha = \mp (1-\alpha^{2m})^{\frac{1}{2m}} +c_2 \]
for a constant $c_2$. It satisfies
\[\alpha^{2m}+(u_{\pm}(\alpha)-c_2)^{2m} = 1. \]
So the resulting surface can be smoothly extended to a Minkowski sphere, which is a parallel translation of the unit sphere $S$. 

\vspace{2mm}

(i-2) The case $0 < c_1 < 1$. In this case, we have $0 < \alpha < \sqrt{1-c_1}$ and we can write
\[u_{\pm}(\alpha) = \pm \int_{\sqrt{1-c_1}}^{\alpha} \ \frac{(\rho^2+c_1)^{\frac{2m-1}{2}}} {\left\{ 1-(\rho^2+c_1)^m \right\}^{\frac{2m-1}{2m}}} \ d\rho+c_3 \]
for a constant $c_3$. Since
\[0 < \frac{2m-1}{2m} < 1, \]
the above integral converges. Set
\[d_1 := -\lim_{\alpha\rightarrow 0} \int_{\sqrt{1-c_1}}^{\alpha} \ \frac{(\rho^2+c_1)^{\frac{2m-1}{2}}} {\left\{ 1-(\rho^2+c_1)^m \right\}^{\frac{2m-1}{2m}}} \ d\rho \ \ (> 0). \]
Then
\[\lim_{\alpha\rightarrow 0} u_{+}(\alpha) = c_3-d_1, \ \ \ \ \lim_{\alpha\rightarrow 0} u_{-}(\alpha) = c_3+d_1, \ \ \ \ \lim_{\alpha\rightarrow \sqrt{1-c_1}} u_{\pm}(\alpha) = c_3. \]
The function $u_{+}(\alpha)$ is an increasing function and
\[\lim_{\alpha\rightarrow \sqrt{1-c_1}} u_{+}'(\alpha) = \infty. \]
The function $u_{-}(\alpha)$ is a decreasing function and
\[\lim_{\alpha\rightarrow \sqrt{1-c_1}} u_{-}'(\alpha) = -\infty. \]
Let $\alpha_{+}(u)$ be the inverse function of $u_{+}(\alpha)$. It is an increasing function on $(c_3-d_1, c_3)$ and
\[\lim_{u\rightarrow c_3-d_1} \alpha_{+}(u) = 0, \ \ \ \ \lim_{u\rightarrow c_3} \alpha_{+}(u) = \sqrt{1-c_1}, \ \ \ \ \lim_{u\rightarrow c_3} \alpha_{+}'(u) = 0. \]
Let $\alpha_{-}(u)$ be the inverse function of $u_{-}(\alpha)$. It is a decreasing function on $(c_3, c_3+d_1)$ and
\[\lim_{u\rightarrow c_3+d_1} \alpha_{-}(u) = 0, \ \ \ \ \lim_{u\rightarrow c_3} \alpha_{-}(u) = \sqrt{1-c_1}, \ \ \ \ \lim_{u\rightarrow c_3} \alpha_{-}'(u) = 0. \]

We define a function $\hat{\alpha}(u)$ on $(c_3-d_1, c_3+d_1)$ by
\[\hat{\alpha}(u) = \left\{ \begin{array}{cc} 
\alpha_{+}(u), & c_3-d_1 < u < c_3 \\ \\
\alpha_{-}(u), & c_3 < u < c_3+d_1 \\ \\
\sqrt{1-c_1}, & u = c_3 
\end{array}. \right. \]
Then $\hat{\alpha}(u)$ is a $C^1$-function on $(c_3-d_1, c_3+d_1)$ such that
\[\hat{\alpha}'(u) = \left\{ \begin{array}{cc} 
\alpha_{+}'(u), & c_3-d_1 < u < c_3 \\ \\
\alpha_{-}'(u), & c_3 < u < c_3+d_1 \\ \\
0, & u = c_3 
\end{array}. \right. \]
For $u \in (c_3-d_1, c_3) \cup (c_3, c_3+d_1)$, $\hat{\alpha}(u)$ satisfies the equation (5.1) for $K = 1$. Then, noting that $m \geq 2$, we can see that
\[\lim_{u\rightarrow c_3} (\hat{\alpha}'(u))^{-\frac{2m-2}{2m-1}} \hat{\alpha}''(u) = -(2m-1)\sqrt{1-c_1} \]
and
\[\lim_{u\rightarrow c_3} \hat{\alpha}''(u) = 0. \]
Thus the function $\hat{\alpha}(u)$ is a $C^2$-function on $(c_3-d_1, c_3+d_1)$. Also by (3.5) and (3.6), we find that the Birkhoff-Gauss map can be $C^1$-extended for $u \in (c_3-d_1, c_3+d_1)$. 

On the other hand, we have
\[\lim_{u\rightarrow c_3-d_1} \hat{\alpha}(u) = 0, \ \ \ \ \lim_{u\rightarrow c_3+d_1} \hat{\alpha}(u) = 0 \]
and
\[\lim_{u\rightarrow c_3-d_1} \hat{\alpha}'(u) = \frac{(1-c_{1}^{m})^{\frac{2m-1}{2m}}} {c_{1}^{\frac{2m-1}{2}}}, \ \ \ \ \lim_{u\rightarrow c_3+d_1} \hat{\alpha}'(u) = -\frac{(1-c_{1}^{m})^{\frac{2m-1}{2m}}} {c_{1}^{\frac{2m-1}{2}}}. \]
So the surface has singularities at $(0,0,c_3-d_1)$ and $(0,0,c_3+d_1)$. 

\vspace{2mm}

(i-3) The case $c_1 < 0$. In this case, we have $\sqrt{-c_1} < \alpha < \sqrt{1-c_1}$ and
\[u_{\pm}(\alpha) = \pm \int_{\sqrt{1-c_1}}^{\alpha} \ \frac{(\rho^2+c_1)^{\frac{2m-1}{2}}} {\left\{ 1-(\rho^2+c_1)^m \right\}^{\frac{2m-1}{2m}}} \ d\rho+c_4 \]
for a constant $c_4$. As in the case (i-2), we can see that the graphs of $u_{+}(\alpha)$ and $u_{-}(\alpha)$ are connected smoothly at $\alpha = \sqrt{1-c_1}$. But the surface has singularities at points where $\alpha = \sqrt{-c_1}$. 

\vspace{2mm}

(ii) The case $K = -1$. We have $c_1 > 0$ and 
\[u_{\pm}(\alpha) = \pm \int \frac{(c_1-\alpha^2)^{\frac{2m-1}{2}}} {\left\{ 1-(c_1-\alpha^2)^m \right\}^{\frac{2m-1}{2m}}} d\alpha. \]

(ii-1) The case $0 < c_1 \leq 1$. In this case, we have $0 < \alpha < \sqrt{c_1}$ and
\[u_{\pm}(\alpha) = \pm \int_{\sqrt{c_1}}^{\alpha} \frac{(c_1-\rho^2)^{\frac{2m-1}{2}}} {\left\{ 1-(c_1-\rho^2)^m \right\}^{\frac{2m-1}{2m}}} d\rho+c_5 \]
for a constant $c_5$. Then
\[\lim_{\alpha\rightarrow \sqrt{c_1}} u_{\pm}(\alpha) = c_5, \ \ \ \ \lim_{\alpha\rightarrow \sqrt{c_1}} u_{\pm}'(\alpha) = 0. \]

(ii-1-1) When $c_1 = 1$, since
\[1 < \frac{2m-1}{m} < 2, \]
we have
\[\lim_{\alpha\rightarrow 0} u_{\pm}(\alpha) = \mp \infty. \]
The corresponding surface has singularities at points where $\alpha = \sqrt{c_1}$, and it can be seen as a generalization of the pseudo-sphere in the Euclidean $3$-space. 

(ii-1-2) When $0 < c_1 < 1$, we have
\[\lim_{\alpha\rightarrow 0} u_{\pm}(\alpha) = c_5 \mp d_2 \]
where
\[d_2 := -\lim_{\alpha\rightarrow 0} \int_{\sqrt{c_1}}^{\alpha} \frac{(c_1-\rho^2)^{\frac{2m-1}{2}}} {\left\{ 1-(c_1-\rho^2)^m \right\}^{\frac{2m-1}{2m}}} d\rho \ \ (> 0), \]
and
\[\lim_{\alpha\rightarrow 0} u_{\pm}'(\alpha) = \pm \frac{c_{1}^{\frac{2m-1}{2}}} {(1-c_1^m)^{\frac{2m-1}{2m}}}. \]
So the surface has singularities at points where $\alpha = \sqrt{c_1}$ and $\alpha = 0$. 

\vspace{3mm}

(ii-2) The case $c_1 > 1$. In this case, we have $\sqrt{c_1-1} < \alpha < \sqrt{c_1}$ and
\[u_{\pm}(\alpha) = \pm \int_{\sqrt{c_1-1}}^{\alpha} \frac{(c_1-\rho^2)^{\frac{2m-1}{2}}} {\left\{ 1-(c_1-\rho^2)^m \right\}^{\frac{2m-1}{2m}}} d\rho+c_6 \]
for a constant $c_6$. By the discussion as before, the graphs of $u_{+}$ and $u_{-}$ can be $C^2$-connected at $\alpha = \sqrt{c_1-1}$. But the surface has singularities at points where $\alpha = \sqrt{c_1}$.

\section{Non-zero constant mean curvature}

Let $({\mathbb R}^3, \|\cdot\|)$ be the $3$-dimensional normed space as in Section 3. Let $M$ be a rotational surface in $({\mathbb R}^3, \|\cdot\|)$ parametrized by
\[f(u, v) = (\alpha(u)\cos{v}, \alpha(u)\sin{v}, u) \]
where $\alpha > 0$ and $\alpha' \neq 0$. 

By (3.8), if $H$ is a non-zero constant, then
\begin{eqnarray}
\frac{1}{2m-1} \alpha \left( \left( \alpha' \right)^{\frac{2m}{2m-1}}+1 \right)^{-\frac{2m+1}{2m}} (\alpha')^{-\frac{2m-2}{2m-1}} \alpha'' -\left( \left( \alpha' \right)^{\frac{2m}{2m-1}}+1 \right)^{-\frac{1}{2m}} \nonumber
\end{eqnarray}
\begin{eqnarray}
\hspace{1cm} = 2H\alpha. 
\end{eqnarray}
Multiplying by $-\alpha'$ we have
\[-\frac{1}{2m-1} \alpha \left( \left( \alpha' \right)^{\frac{2m}{2m-1}}+1 \right)^{-\frac{2m+1}{2m}} (\alpha')^{\frac{1}{2m-1}} \alpha'' +\alpha' \left( \left( \alpha' \right)^{\frac{2m}{2m-1}}+1 \right)^{-\frac{1}{2m}} = -H(\alpha^2)'. \]
Integrating it we can get
\[\alpha \left( \left( \alpha' \right)^{\frac{2m}{2m-1}}+1 \right)^{-\frac{1}{2m}} = c_1-H\alpha^2 \ \ (> 0) \]
for a constant $c_1$. Then
\[\frac{d\alpha}{du} = \pm \frac{\left\{ \alpha^{2m}-(c_1-H\alpha^2)^{2m} \right\}^{\frac{2m-1}{2m}}} {(c_1-H\alpha^2)^{2m-1}}, \]
and we get the following. 

\begin{thm}
A rotational surface in $({\mathbb R}^3, \|\cdot\|)$ given by
\[\bar{f}(\alpha, v) = (\alpha \cos v, \alpha \sin v, u(\alpha)) \]
where $\alpha > 0$ has non-zero constant Minkowski mean curvature $H$ if and only if
\[u(\alpha) = \pm \int \frac{(c_1-H\alpha^2)^{2m-1}} {\left\{ \alpha^{2m}-(c_1-H\alpha^2)^{2m} \right\}^{\frac{2m-1}{2m}}} d\alpha \]
for a constant $c_1$. 
\end{thm}

Set
\[u_{\pm}(\alpha) := \pm \int \frac{(c_1-H\alpha^2)^{2m-1}} {\left\{ \alpha^{2m}-(c_1-H\alpha^2)^{2m} \right\}^{\frac{2m-1}{2m}}} d\alpha. \]
To study the behavior of the graph of $u_{\pm}(\alpha)$, it suffices to consider the case where $H = \pm1$. The signature of $H$ changes if the orientation of the parametrization changes. So we treat the both cases $H = 1$ and $H = -1$. 

\vspace{2mm}

(i) The case $H = 1$. In this case, we have $c_1 > 0$, 
\[b_1 := \frac{\sqrt{1+4c_1}-1}{2} < \alpha < \sqrt{c_1} \]
and
\[u_{\pm}(\alpha) = \pm \int_{\sqrt{c_1}}^{\alpha} \frac{(c_1-\rho^2)^{2m-1}} {\left\{ \rho^{2m}-(c_1-\rho^2)^{2m} \right\}^{\frac{2m-1}{2m}}} d\rho+c_{2}^{\pm} \]
for a constant $c_{2}^{\pm}$. This integral converges as $\alpha$ tends to $\sqrt{c_1}$, and since $0 < (2m-1)/2m < 1$, it converges also as $\alpha$ tends to $b_1$. Set
\[d_1 := -\lim_{\alpha\rightarrow b_1} \int_{\sqrt{c_1}}^{\alpha} \frac{(c_1-\rho^2)^{2m-1}} {\left\{ \rho^{2m}-(c_1-\rho^2)^{2m} \right\}^{\frac{2m-1}{2m}}} d\rho \ \ (> 0). \]
Then
\[\lim_{\alpha\rightarrow b_1} u_{+}(\alpha) = c_{2}^{+}-d_1, \ \ \ \ \lim_{\alpha\rightarrow b_1} u_{-}(\alpha) = c_{2}^{-}+d_1, \ \ \ \ \lim_{\alpha\rightarrow \sqrt{c_1}} u_{\pm}(\alpha) = c_{2}^{\pm}. \]
The function $u_{+}(\alpha)$ is an increasing function and
\[\lim_{\alpha\rightarrow b_1} u_{+}'(\alpha) = \infty, \ \ \ \ \lim_{\alpha\rightarrow \sqrt{c_1}} u_{+}'(\alpha) = 0. \]
The function $u_{-}(\alpha)$ is a decreasing function and
\[\lim_{\alpha\rightarrow b_1} u_{-}'(\alpha) = -\infty, \ \ \ \ \lim_{\alpha\rightarrow \sqrt{c_1}} u_{-}'(\alpha) = 0. \]

(ii) The case $H = -1$. We have
\[u_{\pm}(\alpha) = \pm \int \frac{(c_3+\alpha^2)^{2m-1}} {\left\{ \alpha^{2m}-(c_3+\alpha^2)^{2m} \right\}^{\frac{2m-1}{2m}}} d\alpha \]
for a constant $c_3$. Here we use $c_3$ instead of $c_1$ because we will later choose $c_3$ different from $c_1$. 

\vspace{2mm}

(ii-1) The case $c_3 = 0$. Then
\[u_{\pm}(\alpha) = \pm \int \frac{\alpha^{2m-1}} {(1-\alpha^{2m})^{\frac{2m-1}{2m}}} d\alpha = \mp(1-\alpha^{2m})^{\frac{1}{2m}}+c_4 \]
for a constant $c_4$. It satisfies
\[\alpha^{2m}+(u_{\pm}(\alpha)-c_4)^{2m} = 1. \]
So the surface can be smoothly extended to a Minkowski sphere, which is a parallel translation of the unit sphere $S$. 

\vspace{2mm}

(ii-2) The case $c_3 > 0$. In this case, we have $0 < c_3 < 1/4$, 
\[b_2 := \frac{1-\sqrt{1-4c_3}}{2} < \alpha < \frac{1+\sqrt{1-4c_3}}{2} =: b_3 \]
and
\[u_{\pm}(\alpha) = \pm \int_{b_2}^{\alpha} \frac{(c_3+\rho^2)^{2m-1}} {\left\{ \rho^{2m}-(c_3+\rho^2)^{2m} \right\}^{\frac{2m-1}{2m}}} d\rho+c_5 \]
for a constant $c_5$. This integral converges as $\alpha$ tends to $b_2$ and $b_3$. Set
\[d_2 := \lim_{\alpha\rightarrow b_3} \int_{b_2}^{\alpha} \frac{(c_3+\rho^2)^{2m-1}} {\left\{ \rho^{2m}-(c_3+\rho^2)^{2m} \right\}^{\frac{2m-1}{2m}}} d\rho. \]
Then
\[\lim_{\alpha\rightarrow b_2} u_{\pm}(\alpha) = c_5, \ \ \ \ \lim_{\alpha\rightarrow b_3} u_{+}(\alpha) = c_5+d_2, \ \ \ \ \lim_{\alpha\rightarrow b_3} u_{-}(\alpha) = c_5-d_2. \]
The function $u_{+}(\alpha)$ is an increasing function and
\[\lim_{\alpha\rightarrow b_2} u_{+}'(\alpha) = \infty, \ \ \ \ \lim_{\alpha\rightarrow b_3} u_{+}'(\alpha) = \infty. \]
The function $u_{-}(\alpha)$ is a decreasing function and
\[\lim_{\alpha\rightarrow b_2} u_{-}'(\alpha) = -\infty, \ \ \ \ \lim_{\alpha\rightarrow b_3} u_{-}'(\alpha) = -\infty. \]
Let $\alpha_{+}(u)$ be the inverse function of $u_{+}(\alpha)$. It is increasing on $(c_5, c_5+d_2)$ and
\[\lim_{u\rightarrow c_5} \alpha_{+}(u) = b_2, \ \ \ \ \lim_{u\rightarrow c_5+d_2} \alpha_{+}(u) = b_3, \ \ \ \ \lim_{u\rightarrow c_5} \alpha_{+}'(u) = \lim_{u\rightarrow c_5+d_2} \alpha_{+}'(u) = 0. \]
Let $\alpha_{-}(u)$ be the inverse function of $u_{-}(\alpha)$. It is decreasing on $(c_5-d_2, c_5)$ and
\[\lim_{u\rightarrow c_5} \alpha_{-}(u) = b_2, \ \ \ \ \lim_{u\rightarrow c_5-d_2} \alpha_{-}(u) = b_3, \ \ \ \ \lim_{u\rightarrow c_5} \alpha_{-}'(u) = \lim_{u\rightarrow c_5-d_2} \alpha_{-}'(u) = 0. \]

We define a function $\hat{\alpha}(u)$ on $[c_5-d_2, c_5+d_2]$ by
\[\hat{\alpha}(u) = \left\{ \begin{array}{cc} 
\alpha_{+}(u), & c_5 < u < c_5+d_2 \\ \\
\alpha_{-}(u), & c_5-d_2 < u < c_5 \\ \\
b_2, & u = c_5 \\ \\
b_3, & u = c_5 \pm d_2 
\end{array}. \right. \]
Then it is a $C^1$-function on $[c_5-d_2, c_5+d_2]$ such that
\[\hat{\alpha}'(u) = \left\{ \begin{array}{cc} 
\alpha_{+}'(u), & c_5 < u < c_5+d_2 \\ \\
\alpha_{-}'(u), & c_5-d_2 < u < c_5 \\ \\
0, & u = c_5, \ c_5 \pm d_2 
\end{array}. \right. \]
For $u \in (c_5-d_2, c_5) \cup (c_5, c_5+d_2)$, $\hat{\alpha}(u)$ satisfies the equation (6.1) for $H = -1$. Then, noting that $m \geq 2$, we can see that
\[\lim_{u\rightarrow c_5} (\hat{\alpha}'(u))^{-\frac{2m-2}{2m-1}} \hat{\alpha}''(u) = \frac{(2m-1)(1-2b_2)}{b_2}, \]
\[\lim_{u\rightarrow c_5 \pm d_2} (\hat{\alpha}'(u))^{-\frac{2m-2}{2m-1}} \hat{\alpha}''(u) = \frac{(2m-1)(1-2b_3)}{b_3} \]
and
\[\lim_{u\rightarrow c_5} \hat{\alpha}''(u) = \lim_{u\rightarrow c_5 \pm d_2} \hat{\alpha}''(u) = 0. \]
So the function $\hat{\alpha}(u)$ is a $C^2$-function on $[c_5-d_2, c_5+d_2]$. By (3.5) and (3.6), the Birkhoff-Gauss map can be $C^1$-extended for $u \in [c_5-d_2, c_5+d_2]$. 

We note that $\hat{\alpha}(u)$ has the same derivatives at the end points $u = c_5-d_2$ and $u = c_5+d_2$. Thus we can extend $\hat{\alpha}(u)$ periodically as a $C^2$-function on ${\mathbb R}$ as follows: 
\[\alpha^{\ast}(u+2k d_2) := \hat{\alpha}(u), \ \ \ \ u \in [c_5-d_2, c_5+d_2], \ \ k \in {\mathbb Z}, \]
and we get the following. 

\begin{thm}
Under the notation above, the rotational surface in $({\mathbb R}^3, \|\cdot\|)$ parametrized by
\[f^{\ast}(u, v) = (\alpha^{\ast}(u)\cos v, \alpha^{\ast}(u)\sin v, u), \ \ \ \ (u, v) \in {\mathbb R} \times [0, 2\pi] \]
has constant Minkowski mean curvature $-1$. 
\end{thm}

Remark. The surface in Theorem 6.2 can be seen as a generalization of the unduloid (\cite{D}). 

\vspace{3mm}

(ii-3) The case $c_3 < 0$. In this case we have
\[\sqrt{-c_3} < \alpha < \frac{1+\sqrt{1-4c_3}}{2} =: b_4 \]
and
\[u_{\pm}(\alpha) = \pm \int_{b_4}^{\alpha} \frac{(c_3+\rho^2)^{2m-1}} {\left\{ \rho^{2m}-(c_3+\rho^2)^{2m} \right\}^{\frac{2m-1}{2m}}} d\rho+c_6 \]
for a constant $c_6$. This integral converges as $\alpha$ tends to $\sqrt{-c_3}$ and $b_4$. Set
\[d_3 := -\lim_{\alpha\rightarrow \sqrt{-c_3}} \int_{b_4}^{\alpha} \frac{(c_3+\rho^2)^{2m-1}} {\left\{ \rho^{2m}-(c_3+\rho^2)^{2m} \right\}^{\frac{2m-1}{2m}}} d\rho \ \ (> 0). \]
Then
\[\lim_{\alpha\rightarrow \sqrt{-c_3}} u_{+}(\alpha) = c_6-d_3, \ \ \ \ \lim_{\alpha\rightarrow \sqrt{-c_3}} u_{-}(\alpha) = c_6+d_3, \ \ \ \ \lim_{\alpha\rightarrow b_4} u_{\pm}(\alpha) = c_6. \]
The function $u_{+}(\alpha)$ is an increasing function and
\[\lim_{\alpha\rightarrow \sqrt{-c_3}} u_{+}'(\alpha) = 0, \ \ \ \ \lim_{\alpha\rightarrow b_4} u_{+}'(\alpha) = \infty. \]
The function $u_{-}(\alpha)$ is a decreasing function and
\[\lim_{\alpha\rightarrow \sqrt{-c_3}} u_{-}'(\alpha) = 0, \ \ \ \ \lim_{\alpha\rightarrow b_4} u_{-}'(\alpha) = -\infty. \]

In the following, we will connect the curves in the cases (i) and (ii-3). For distinguishment, we denote $u_{\pm}(\alpha)$ in the case (i) by $u_{1\pm}(\alpha)$, and $u_{\pm}(\alpha)$ in the case (ii-3) by $u_{2\pm}(\alpha)$. 

We take the graph $G_1$ of $u_{1+}(\alpha)$ for $b_1 < \alpha < \sqrt{c_1}$. Next, choosing $c_3 := -c_1$ and $c_6 := c_{2}^{+}-d_3$, we take the graph $G_2$ of $u_{2-}(\alpha)$ for 
\[\sqrt{-c_3} = \sqrt{c_1} < \alpha < b_4 = \frac{1+\sqrt{1+4c_1}}{2}. \]
Then $G_1$ and $G_2$ are $C^1$-connected at $(\alpha, u) = (\sqrt{c_1}, c_{2}^{+})$. 

Next we take the graph $G_3$ of $u_{2+}(\alpha)$ for $\sqrt{c_1} < \alpha < b_4$. Then $G_2$ and $G_3$ are $C^1$-connected at $(\alpha, u) = (b_4, c_{2}^{+}-d_3)$. 

Finally, letting $c_{2}^{-} := c_{2}^{+}-2d_3$, we take the graph $G_4$ of $u_{1-}(\alpha)$ for $b_1 < \alpha < \sqrt{c_1}$. Then $G_3$ and $G_4$ are $C^1$-connected at $(\alpha, u) = (\sqrt{c_1}, c_{2}^{+}-2d_3)$. Thus we get a $C^1$-curve $\Gamma$ which is constructed by connecting $G_1$, $G_2$, $G_3$ and $G_4$. 

With respect to the parameter $u$, $H = 1$ for the $G_1$ and $G_4$ parts, and $H = -1$ for the $G_2$ and $G_3$ parts. On the other hand, with respect to the parameter $\alpha$, $H = 1$ for the $G_1$ and $G_2$ parts, and $H = -1$ for the $G_3$ and $G_4$ parts. Then, with respect to a parametrization of $\Gamma$ in the order of $G_1$, $G_2$, $G_3$ and $G_4$, we have $H = 1$ for all parts. 

The $C^2$-connectedness of $G_2$ and $G_3$ is shown by the discussion as before. Similarly we can see that $G_1$ and $G_4$ are $C^2$ at $\alpha = b_1$. 

Let us prove the $C^2$-connectedness of $G_1$ and $G_2$. We define a function $\hat{u}(\alpha)$ on $(b_1, b_4)$ by
\[\hat{u}(\alpha) = \left\{ \begin{array}{cc} 
u_{1+}(\alpha), & b_1 < \alpha < \sqrt{c_1} \\ \\
u_{2-}(\alpha), & \sqrt{c_1} < \alpha < b_4 \\ \\
c_{2}^{+}, & \alpha = \sqrt{c_1} 
\end{array}. \right. \]
Then it is a $C^1$-function on $(b_1, b_4)$ such that
\[\hat{u}'(\alpha) = \left\{ \begin{array}{cc} 
u_{1+}'(\alpha), & b_1 < \alpha < \sqrt{c_1} \\ \\
u_{2-}'(\alpha), & \sqrt{c_1} < \alpha < b_4 \\ \\
0, & \alpha = \sqrt{c_1} 
\end{array}. \right. \]
For $\alpha \in (b_1, \sqrt{c_1}) \cup (\sqrt{c_1}, b_4)$, $\hat{u}(\alpha)$ satisfies the equation (3.4) for "$H = 1$", where $\alpha$ is the parameter and $\beta = \hat{u}(\alpha)$. Then, noting that $m \geq 2$, we can see that
\[\lim_{\alpha\rightarrow \sqrt{c_1}} \left( \hat{u}'(\alpha) \right)^{-\frac{2m-2}{2m-1}} \hat{u}''(\alpha) = -2(2m-1) \]
and
\[\lim_{\alpha\rightarrow \sqrt{c_1}}  \hat{u}''(\alpha) = 0. \]
So the function $\hat{u}(\alpha)$ is a $C^2$-function on $(b_1, b_4)$. By (3.1) and (3.2), the Birkhoff-Gauss map can be $C^1$-extended for $\alpha \in (b_1, b_4)$. Thus the $C^2$-connectedness of $G_1$ and $G_2$ is proved. The $C^2$-connectedness of $G_3$ and $G_4$ is proved similarly. 

Now we have obtained a $C^2$-curve $\Gamma$ which is constructed by connecting $G_1$, $G_2$, $G_3$ and $G_4$. The curve $\Gamma$ has the same derivatives at the end points. Then, as in the case (ii-2), we can extend it periodically as a $C^2$-curve $\Gamma^{\ast}$, which can be parametrized as $(\alpha^{\ast}(t), \beta^{\ast}(t))$ for $t \in {\mathbb R}$. 

\begin{thm}
Under the notation above, the rotational surface in $({\mathbb R}^3, \|\cdot\|)$ parametrized by
\[f^{\ast}(t, v) = (\alpha^{\ast}(t)\cos v, \alpha^{\ast}(t)\sin v, \beta^{\ast}(t)), \ \ \ \ (t, v) \in {\mathbb R} \times [0, 2\pi] \]
has constant Minkowski mean curvature $1$. 
\end{thm}

Remark. (i) The surface in Theorem 6.3 can be seen as a generalization of the nodoid (\cite{D}). 

(ii) By "Mathematica" we know that: (a) When $c_1 = 2$ and $m = 2$, $d_1 = 0.34459...$ and $d_3 = 0.65540...$, (b) When $c_1 = 2$ and $m = 3$, $d_1 = 0.33886...$ and $d_3 = 0.66113...$, and (c) When $c_1 = 6$ and $m = 2$, $d_1 = 0.40710...$ and $d_3 = 0.59289...$. Thus the curve $\Gamma$ is not closed in those cases.

\vspace{3mm}

{\small 
Makoto SAKAKI

Graduate School of Science and Technology, Hirosaki University

Hirosaki 036-8561, Japan

E-mail: sakaki@hirosaki-u.ac.jp
}

\end{document}